\documentclass[a4paper,12pt,intlimits,oneside]{amsart}
\usepackage{amsmath}
\usepackage{amsthm}
\usepackage{latexsym}
\usepackage{amssymb}
\usepackage{mathrsfs}
\usepackage{xcolor}
\usepackage{graphicx}
\usepackage[colorlinks=true]{hyperref}
\numberwithin{figure}{section}
\def\R{{\mathbb R}}
\def\C{{\mathbb C}}
\def\D{{\mathbb D}}
\def\T{{\mathbb T}}

\def\la{\langle}
\def\ra{\rangle}
\def\s{\vskip 0.25cm\noindent}
\def\e{\varepsilon}

\def\build#1_#2^#3{\mathrel{
\mathop{\kern 0pt#1}\limits_{#2}^{#3}}}
\def\td_#1,#2{\mathrel{\mathop{\build\longrightarrow_{#1\rightarrow #2}^{}}}}
\DeclareFontFamily{U}{MnSymbolC}{}
\DeclareSymbolFont{MnSyC}{U}{MnSymbolC}{m}{n}
\DeclareFontShape{U}{MnSymbolC}{m}{n}{
    <-6>  MnSymbolC5
   <6-7>  MnSymbolC6
   <7-8>  MnSymbolC7
   <8-9>  MnSymbolC8
   <9-10> MnSymbolC9
  <10-12> MnSymbolC10
  <12->   MnSymbolC12}{}
\DeclareMathSymbol{\intprod}{\mathbin}{MnSyC}{'270}
\newtheorem{theorem}{Theorem}

\newtheorem{lemma}{Lemma}
\newtheorem{remark}{Remark}

\begin{document}
\title[An explicit formula for the Benjamin--Ono equation]{An explicit formula for the Benjamin--Ono equation}
\author[P. G\'erard]{Patrick G\'erard}
\address{Universit\'e Paris--Saclay, Laboratoire de Math\'ematiques d'Orsay,  CNRS,   UMR 8628,  91405 Orsay, France} \email{{\tt patrick.gerard@universite-paris-saclay.fr}}


\date{December 5,  2022}

\begin{abstract}
We establish an explicit formula for the general solution of the Benjamin--Ono equation on the torus and on the line.
\end{abstract}


\thanks{I would like to warmly thank my colleagues
M. Bellassoued, N. Burq,  B. Dehman and C. Zuily 
for inviting me to give a minicourse in the Tunisian--French conference
of October 2022 in Hammamet, which encouraged me to provide a joint presentation of the cubic Szeg\H{o} equation and of the Benjamin--Ono equation,
and led me to the explicit formulae established here. I am grateful to Sandrine Grellier, Enno Lenzmann, Peter Perry and Peter Topalov   for fruitful discussions about these formulae. }

\maketitle

\tableofcontents

\medskip

\section{Introduction}
\subsection{The Benjamin--Ono equation}
The Benjamin--Ono equation was introduced by Benjamin \cite{B} (see also Davis--Acrivos \cite{DA}, Ono \cite{O}) to model long, one-way internal gravity waves 
in a two-layer fluid.
It  reads 
\begin{equation} \label{BO}
\partial_tu=\partial_x(|D|u-u^2)\ .
\end{equation}
Here $u=u(t,x)$ denotes a real valued function. There is a vast literature about  this equation, and we refer to the book by Klein and Saut \cite{KS} for a recent survey. We consider both the case of periodic boundary conditions $u(t,x+2\pi )=u(t,x)$, which we refer as $x\in \T$, and the case where $u(t,x)$ cancels as $x$ tends to $\pm \infty $, which we refer as $x\in \R$.  In both cases,  we will restrict ourselves to sufficiently smooth solutions, which can be proved to exist globally by a combination of standard quasilinear scheme and appropriate conservation laws. More precisely, if we denote by 
 $H^2_r$ the Sobolev space of real valued functions with two derivatives in $L^2$, one can prove the following result.
\begin{theorem}[Saut, 1979 \cite{Saut}]\label{JCS} 
For every $u_0\in H^2_r$, there exists a unique  solution $u\in C(\R ,H^2_r)$ of \eqref{BO} with $u(0)=u_0$.
\end{theorem}
Our goal in this paper is to provide an explicit formula of the solution $u(t)$ in terms of the initial datum $u_0$.
\s
For this, we need to  introduce the Lax pair structure for \eqref{BO}. 
\subsection{The Lax pair}
On $\T $ or $\R $, we denote by $L^2_+$  the Hardy space corresponding to $L^2$ functions having a Fourier transform supported in the domain $\xi \geq 0$. Recall that both Hardy spaces identify to some spaces of holomorphic functions. 
The space $L^2_+(\T )$ identifies to holomorphic functions $f$ on the unit disc $\D :=\{ z\in \C : |z|<1\}$ such that 
$$\sup_{r<1} \int_0^{2\pi}\left |f(r{\mathrm e}^{ix})\right |^2\, dx <+\infty ,$$
while $L^2_+(\R )$ identifies to holomorphic functions on the upper half plane $\C_+:=\{ z\in \C : \mathrm{Im}(z)>0\} $ such that
$$\sup_{y>0}\int_\R |f(x+iy)|^2\, dy <+\infty \ .$$
We denote by $\Pi $ the orthogonal projector from $L^2$ onto $L^2_+$. 
Remarkable operators on $L^2_+$ are Toeplitz operators, associated to functions $b\in L^\infty $ by the formula
$$\forall f\in L^2_+\ ,\ T_b f=\Pi (bf)\ .$$
For every $u\in L^2_r$, we denote by $L_u$ the semi--bounded selfadjoint operator defined on $L^2_+$ by 
$$\mathrm {Dom}(L_u)=H^1_+:=H^1\cap L^2_+\ ,\ L_uf=Df-T_uf\ ,\ D:=\frac{1}{i}\frac{d}{dx}\ .$$
We also consider, for $u\in H^2_r$, the bounded antiselfadjoint operator defined by
$$B_u=i(T_{|D|u}-T_u^2)\ .$$
Then one can check the following identity (see e.g. \cite{FA}, \cite{Wu}, \cite{GK}, \cite{G}).
\begin{theorem}\label{Laxpair} 
Under the conditions of Theorem \ref{JCS}, we have 
$$\forall t\in \R \ ,\ \frac{d}{dt}L_{u(t)}=[B_{u(t)}, L_{u(t)}]\ .$$
\end{theorem}
\subsection{The explicit formula on the torus}
Let us mention some more properties of  the Hardy space on the torus. The Hardy space $L^2_+(\T )$ is equipped with the shift operator $S:=T_{\mathrm{e}^{ix}}$ and with its adjoint $S^*=T_{\mathrm{e}^{-ix}}$. 
With this notation, our main result on the torus reads as follows.
\begin{theorem}\label{torus}
Let $u\in C(\R, H^2_r(\T ))$ be the solution of the Benjamin--Ono equation on the torus $\T $ with $u(0)=u_0$. \\
Then $u(t)=\Pi u (t)+\overline {\Pi u}(t) -\la u_0\vert 1\ra ,$ with
$$\forall z\in \D \ ,\ \Pi u(t,z)=\left \la  (\mathrm {Id}-z{\mathrm e}^{it}{\mathrm e}^{2itL_{u_0}}S^*)^{-1} \Pi u_0\vert 1  \right \ra \ .$$
\end{theorem}
\begin{remark}\label{remtor}
The above formula can be equivalently stated as a characterization of Fourier coefficients of the solution $u$,
\begin{eqnarray*}
\forall k\geq 0\ ,\ \hat u(t,k)&=&\left \la  ({\mathrm e}^{it}{\mathrm e}^{2itL_{u_0}}S^*)^k \Pi u_0\vert  1\right \ra \\
&=&\left \la \Pi u_0\vert (S{\mathrm e}^{-it}{\mathrm e}^{-2itL_{u_0}})^k\, 1\right \ra  \ .
\end{eqnarray*}
Under this form, it extends to much more singular data, for which the flow of the Benjamin--Ono has been proved to extend continuously. 
According to \cite{GKT}, this is the case if $u_0$ belongs to the Sobolev space $H^s_r(\T )$ for every $s>-1/2$. Indeed, in this case, $L_{u_0}$ is selfadjoint, semibounded, 
and the domain of the square root of $L_{u_0}+K{\rm Id}$, for $K$ big enough, is the space $H^{1/2}_+(\T )$. Consequently, the operator 
$S{\mathrm e}^{-it}{\mathrm e}^{-2itL_{u_0}}$ acts on $H^{1/2}_+(\T )$, so that the  inner product in the second line above is well defined.
\end{remark}
\subsection{The explicit formula on the line}
On $L^2_+(\R )$, the shift operator $S$ must be replaced by the Lax--Beurling semi--group $(S(\eta ))_{\eta \geq 0}$ of isometries defined as
$$\forall f\in L^2_+(\R )\ ,\ S(\eta )f (x)={\mathrm e}^{i\eta x}f(x)\ .$$
We denote by  $G$  the adjoint of the operator of multiplication by $x$ on $L^2_+(\R )$. Notice that $-iG$ is the infinitesimal generator of the adjoint semi--group of contractions $(S(\eta )^*)_{\eta \geq 0}$, so that
$$\forall \eta \geq 0\ ,\ S(\eta )^*=\mathrm{e}^{-i\eta G}\ .$$
It is easy to check that 
the domain of $G$ consists of those functions $f\in L^2_+(\R )$ such that the restriction of $\hat f$ to the half--line $(0,+\infty )$ belongs to the Sobolev space $H^1(0,+\infty )$, and that 
$$\widehat {Gf}(\xi )=i\frac{d}{d\xi}[\hat f(\xi )]\, {\bf 1}_{\xi >0}\ .$$
As a consequence, for every $f\in \mathrm{Dom}(G)$, one may define
$$I_+(f):=\hat f(0^+)\ .$$
This definition can be extended to any $f\in L^2_+$ such that the restriction of $\hat f$ to some interval $(0,\delta )$ belongs to the Sobolev space $H^1(0,\delta )$ for some $\delta >0$, and we shall use it as well.
\\
With this notation, our main result on the line reads as follows.
\begin{theorem}\label{line}
Let $u\in C(\R, H^2_r(\R ))$ be the solution of the Benjamin--Ono equation on the line $\R $ with $u(0)=u_0$. \\
Then $u(t)=\Pi u (t)+\overline {\Pi u}(t),$ with
$$\forall z\in \C_+ \ ,\ \Pi u(t,z)=\frac{1}{2i\pi}I_+[(G-2tL_{u_0}-z\mathrm {Id})^{-1} \Pi u_0]\ .$$
\end{theorem}
Notice that, in the above formula, the function $$f_{z,t}:=(G-2tL_{u_0}-z\mathrm {Id})^{-1} \Pi u_0$$ belongs to the domain of $G-2tL_{u_0}$ --- see the end of section \ref{proof:line} for more detail ---, and therefore its Fourier transform satisfies $\hat f\in H^1(0,\delta )$ for every finite $\delta >0$, hence one can define $I_+(f_{z,t})$.
\begin{remark}\label{remline} At this time, the wellposedness theory for \eqref{BO} on the line is slightly less advanced than on the torus, see 
\cite{MP} for a detailed account of this, with extension of the flow map to $L^2_r(\R )$. However, one can easily prove -- see section 3 below --- that the above formula makes sense for $u_0$ in the space $L^\infty (\R )\cap L^2_r(\R )$.
\end{remark}
\subsection{Organization of the paper} Sections \ref{proof:torus} and \ref{proof:line} are respectively devoted to the proofs of Theorems \ref{torus} and \ref{line}. The main idea is to take advantage of commutation identities between the operators of the shift structure of the Hardy space and the operators $L_u$ and $B_u$ of the Lax pair, in the spirit of what  was done in \cite{GG} for the cubic Szeg\H{o} equation on the torus. At the end of Section \ref{proof:line}, we also provide a short discussion of the meaning of the formula, leading to an extension to initial data in $L^\infty (\R )\cap L_r^2(\R )$. Section \ref{final} briefly draws possible applications and extensions to other equations.
\section{Proof of the explicit formula on the torus}\label{proof:torus}
The proof is based on the following lemma.
\begin{lemma}\label{bracket1}
For every $u\in H^2_r(\T )$, 
$$[S^*,B_u]=i((L_u+\mathrm {Id})^2S^*-S^*L_u^2)\ .$$
\end{lemma}
Let us postpone the proof of Lemma \ref{bracket1} and complete the proof of Theorem \ref{torus}.
Since $u(t)$ is real valued, we have the identity
$$u(t)=\Pi u (t)+\overline {\Pi u}(t) -\la u (t)\vert 1\ra ,$$
and $\la u(t)\vert 1\ra =\la u_0\vert 1\ra $ because of the equation. It remains to identify $\Pi u(t)$ as a holomorphic function on the disc.
For this, we proceed as in \cite{GG}, where a similar formula was established for the cubic Szeg\H{o} equation. We have, for every $z\in \D$,
$$\Pi u(t,z)=\sum_{n=0}^\infty z^n \la \Pi u(t), {\mathrm e}^{inx}\ra =\la (\mathrm {Id}-zS^*)^{-1}\Pi u(t)\vert 1\ra \ .$$
We denote by $U=U(t)$ the solution of the linear initial value problem in $\mathscr L(L^2_+(\T ))$,
$$U'(t)=B_{u(t)}U(t)\ ,\ U(0)=\mathrm {Id}\ .$$
Since $B_{u(t)}$ is anti--selfadjoint, $U(t)$ is unitary, and we can write
\begin{equation}\label{eq:Pi u1}
\Pi u(t,z)=\la (\mathrm {Id}-zU(t)^*S^*U(t))^{-1}U(t)^*\Pi u(t)\vert U(t)^*1\ra \ .
\end{equation}
Let us calculate
\begin{eqnarray*}
\frac{d}{dt}U(t)^*1&=&-U(t)^*B_{u(t)}1=-iU(t)^*[(T_{|D|u(t)}-T_{u(t)}^2)(1)]\\
&=&-iU(t)^*[D\Pi u(t)-T_{u(t)}\Pi u(t)]=-iU(t)^*L_{u(t)}\Pi u(t)=iU(t)^*L_{u(t)}^2(1)\\
&=&iL_{u_0}^2U(t)^*1\ ,
\end{eqnarray*}
from which we conclude
$$U(t)^*1={\mathrm e}^{itL_{u_0}^2}(1)\ .$$
Consequently, 
$$U(t)^*\Pi u(t)=-U(t)^*L_{u(t)}(1)=-L_{u_0}U(t)^*(1)=-L_{u_0}{\mathrm e}^{itL_{u_0}^2}(1)={\mathrm e}^{itL_{u_0}^2}(1)\Pi u_0\ .$$
Finally, using Lemma \ref{bracket1},
\begin{eqnarray*}
\frac{d}{dt}U(t)^*S^*U(t)&=&U(t)^*[S^*,B_{u(t)}]U(t)=U(t)^*[i((L_{u(t)}+\mathrm {Id})^2S^*-S^*L_{u(t)}^2]U(t)\\
&=& i (L_{u_0}+\mathrm {Id})^2U(t)^*S^*U(t)-iU(t)^*S^*U(t)L_{u_0}^2\ .
\end{eqnarray*}
from which we infer
$$U(t)^*S^*U(t)={\mathrm e}^{it(L_{u_0}+\mathrm {Id})^2}S^*{\mathrm e}^{-itL_{u_0}^2}\ .$$
Plugging the previous identites into \eqref{bracket1}, we conclude
\begin{eqnarray*}
\Pi u(t,z)&=&\la (\mathrm {Id}-z{\mathrm e}^{it(L_{u_0}+\mathrm {Id})^2}S^*{\mathrm e}^{-itL_{u_0}^2})^{-1}{\mathrm e}^{itL_{u_0}^2}(1)\Pi u_0\vert {\mathrm e}^{itL_{u_0}^2}(1)\ra \\
&=& \la (\mathrm {Id}-z{\mathrm e}^{-itL_{u_0}^2}{\mathrm e}^{it(L_{u_0}+\mathrm {Id})^2}S^*)^{-1}\Pi u_0\vert 1\ra 
\end{eqnarray*}
which yields the claimed formula.
\qed

Finally, let us prove Lemma \ref{bracket1}. First of all, it easy to check the following commutation  identity with the Toeplitz operators,
$$\forall b\in L^\infty(\T ), [T_b,S^*]=\la \ .\ \vert 1 \ra S^*\Pi b\ \ .$$
We infer, using the adjoint Leibniz formula $S^*D=DS^*+S^*$,
\begin{eqnarray*}
[S^*,B_u]&=&i([S^*,T_{|D|u}]-T_u[S^*,T_u]-[S^*,T_u]T_u)\\
&=&i(\la \ .\ \vert 1 \ra S^*D\Pi u -T_u\la \ .\ \vert 1 \ra S^*\Pi u-\la \ .\ \vert 1 \ra S^*\Pi uT_u )\\
&=&i(\la \ .\ \vert 1 \ra (DS^*\Pi u-T_uS^*\Pi u+S^*\Pi u)-\la \ .\ \vert T_u1 \ra S^*\Pi u)\\
&=&i(\la \ .\ \vert 1 \ra (L_uS^*\Pi u+S^*\Pi u)+\la \ .\ \vert L_u1 \ra S^*\Pi u)\\
&=&i((L_u+\mathrm {Id})\la \ .\ \vert 1 \ra S^*\Pi u+(\la \ .\ \vert 1 \ra S^*\Pi u)L_u)\\
&=&i((L_u+\mathrm {Id})( (L_u+\mathrm {Id})S^*-S^*L_u)+((L_u+\mathrm {Id})S^*-S^*L_u)L_u)\\
&=&i((L_u+\mathrm {Id})^2S^*-S^*L_u^2)\ .
\end{eqnarray*}
The proof of Theorem \ref{torus} is complete.\qed

\section{Proof of the explicit formula on the line}\label{proof:line}
 We start with the inverse Fourier transform for every $f\in L^2_+(\R )$, which we can write in the upper--half plane, as an absolutely convergent integral,
$$\forall z\in \C_+\ ,\ f(z)=\frac{1}{2\pi}\int_0^\infty  \mathrm{e}^{iz\xi}\hat f(\xi )\, d\xi \ ,$$
while, in view of the Plancherel theorem,  we have, in $L^2(0,+\infty )$,
$$\hat f(\xi )=\lim_{\e \to 0} \int_\R \mathrm{e}^{-ix\xi}\frac{f(x)}{1+i\e x}\, dx=\lim_{\e \to 0} \la S(\xi)^*f\vert \chi_\e \ra \ ,$$
where $\chi_\e $ denotes the following function in $L^2_+(\R )$,
$$\chi_\e (x):=\frac{1}{1-i\e x}\ .$$
Plugging the second formula into the first one, we infer
\begin{eqnarray*}
f(z)&=&\lim_{\e \to 0}\frac{1}{2\pi}\int_0^\infty  \mathrm{e}^{iz\xi}\la S(\xi)^*f\vert \chi_\e \ra d\xi \\
&=&\lim_{\e \to 0}\frac{1}{2\pi}\int_0^\infty   \mathrm{e}^{iz\xi}\la \mathrm{e}^{-i\xi G}f\vert \chi_\e \ra d\xi \\
&=&\lim_{\e \to 0}\frac{1}{2i\pi} \la (G-z\mathrm{Id})^{-1}f\vert \chi_\e\ra \\
&=&\frac{1}{2i\pi} I_+[(G-z\mathrm{Id})^{-1}f]\ .
\end{eqnarray*}
Since $u(t)$ is real valued, we can write $u(t)=\Pi u(t)+\overline{\Pi u(t)}$, and it remains to characterize $\Pi u(t,z)$ for $z\in \C_+$.
We are going to proceed as in the previous paragraph, using the family $U(t)$ of unitary operators defined by the linear initial value problem in $\mathscr L(L^2_+(\R ))$,
$$U'(t)=B_{u(t)}U(t)\ ,\ U(0)=\mathrm {Id}\ .$$
For every $z\in \C_+$, we have
\begin{eqnarray*}
\Pi u(t,z)&=&\lim_{\e \to 0}\frac{1}{2i\pi} \la U(t)^*(G-z\mathrm{Id})^{-1}\Pi u(t)\vert U(t)^*\chi_\e\ra \\
&=& \lim_{\e \to 0}\frac{1}{2i\pi} \la (U(t)^*GU(t)-z\mathrm{Id})^{-1}U(t)^*\Pi u(t)\vert U(t)^*\chi _\e \ra \ .
\end{eqnarray*}
We have the following lemma.
\begin{lemma}\label{bracket2}
For every $u\in H^2_r(\T )$, 
$$[G,B_u]=-2L_u+i[L_u^2,G]\ .$$
\end{lemma}
Let us postpone the proof of Lemma \ref{bracket2} and complete the proof of Theorem \ref{line}. 
We calculate
\begin{eqnarray*}
\frac{d}{dt}U(t)^*GU(t)&=&U(t)^*[G,B_{u(t)}]U(t)\\
&=&U(t)^*(-2L_{u(t)}+i[L_{u(t)}^2,G])U(t)\\
&=&-2L_{u_0}+i[L_{u_0}^2,U(t)^*GU(t)]\ .
\end{eqnarray*}
Integrating this ODE, we get
$$U(t)^*GU(t)=-2tL_{u_0}+\mathrm{e}^{itL_{u_0}^2}G\mathrm{e}^{-itL_{u_0}^2}\ .$$
Let us determine the other terms in the inner product. We have (see also \cite{Sun})
$$\frac{d}{dt}U(t)^*\Pi u(t)=U(t)^*(\partial_t\Pi u(t)-B_{u(t)}\Pi u(t))=iU(t)^*L_{u(t)}^2\Pi u(t)=iL_{u_0}^2U(t)^*\Pi u(t)\ ,$$
from which we infer
$$U(t)^*\Pi u(t)=\mathrm{e}^{itL_{u_0}^2}\Pi u_0\ .$$
Finally, we have 
$$\frac{d}{dt}U(t)^*\chi_\e =-U(t)^*B_{u(t)}\chi_\e =-iU(t)^*(T_{|D|u(t)}\chi_\e-T_{u(t)}^2\chi_\e )$$ 
and the right hand side converges in $L^2_+$ to
\begin{eqnarray*}&&-iU(t)^*(D\Pi u(t)-T_{u(t)}\Pi u(t) )=-iU(t)^*L_{u(t)}\Pi u(t)\\
&&=-iL_{u_0}U(t)^*\Pi u(t)=-iL_{u_0}\mathrm{e}^{itL_{u_0}^2}\Pi u_0\\
&&=\lim_{\e \to 0}iL_{u_0}^2\mathrm{e}^{itL_{u_0}^2}\chi_\e \ .
\end{eqnarray*}
By integrating in time, we infer
$$U(t)^*\chi_\e -\mathrm{e}^{itL_{u_0}^2}\chi_\e \rightarrow 0$$
in $L^2_+$. Plugging these informations into the formula which gives $\Pi u(t,z)$, we infer
\begin{eqnarray*}
\Pi u(t,z)&=&\lim_{\e \to 0}\frac{1}{2i\pi} \la \left (\mathrm{e}^{itL_{u_0}^2}G\mathrm{e}^{-itL_{u_0}^2}-2tL_{u_0} -z\mathrm{Id}\right )^{-1}
\mathrm{e}^{itL_{u_0}^2}\Pi u_0\vert \mathrm{e}^{itL_{u_0}^2}\chi_\e\ra \\
&=& \lim_{\e \to 0}\frac{1}{2i\pi} \la (G-2tL_{u_0}-z\mathrm{Id})^{-1}\Pi u_0\vert \chi _\e \ra \\
&=& \frac{1}{2i\pi} I_+[(G-2tL_{u_0}-z\mathrm{Id})^{-1}\Pi u_0]\ .
\end{eqnarray*}
It remains to prove Lemma \ref{bracket2}. We shall appeal to the following elementary identity, whose  proof can be found in \cite{Sun}, \cite{GL}. 
\begin{lemma}\label{bracket0}
For every $f\in \mathrm{Dom}(G)$,  $b\in H^1(\R )$, $T_bf\in \mathrm{Dom}(G)$ and
$$[G,T_b]f= \frac{i}{2\pi}I_+(f)\Pi b\ .$$
\end{lemma}
Using  Lemma \ref{bracket0} and the simple observation that 
$[G,D]=i\mathrm{Id}$, we obtain (see also \cite{Sun}),
$$\forall f\in \mathrm{Dom}(G)\cap \mathrm{Dom}(L_u)\ ,\ [G,L_u]f=if-\frac{i}{2\pi}I_+(f)\Pi u\ .$$
We infer, for $f\in \mathrm{Dom}(G)\cap \mathrm{Dom}(L_u^2)$,
\begin{eqnarray*}
[G,B_u]f&=&i([G,T_{|D|u}]f-T_u[G,T_u]f-[G,T_u]T_uf)\\
&=&\frac{i}{2\pi}(iI_+(f)(D\Pi u-T_u\Pi u)-iI_+(T_uf)\Pi u)\\
&=&\frac{i}{2\pi}(iI_+(f)L_u\Pi u+iI_+(L_uf)\Pi u )\\
&=&i(L_u(if -[G,L_u]f)+iL_uf-[G,L_u]L_uf)\\
&=&-2L_uf+i[L_u^2,G]f\ .
\end{eqnarray*}
The proof of Theorem \ref{line} is complete. 
\qed

Let us conclude this section by discussing the formula of Theorem \ref{line} for more singular data $u_0$. First of all, let us observe that, for every $t\in \R$, the operator 
$$A_t:=-i(G-2tL_{0})$$
is maximally dissipative. Indeed, its expression in the Fourier representation is given by 
$$\widehat{A_tf}(\xi )=\frac{d}{d\xi} \hat f(\xi )+2it\xi \hat f(\xi )\ ,$$
and therefore it is easy to check by explicit calculations that 
$$\mathrm{Dom}(A_t)=\{ f\in L^2_+(\R ): {\rm e}^{it\xi ^2}\hat f\in H^1(0,\infty )\}$$
with 
$$\forall f\in \mathrm{Dom}(A_t)\ ,\ {\rm Re}\la A_tf\vert f\ra \leq 0, $$
and that $A_t+i z{\rm Id} : \mathrm{Dom}(A_t)\to L^2_+(\R )$ is bijective for every $z\in \C_+$. From standard perturbation theory, we infer that, for every bounded antiselfadjoint operator $B$ on $L^2_+(\R )$,  $A_t+B$ is maximally dissipative. In particular, if $u_0\in L^\infty (\R)\cap L^2_r(\R )$, 
$$-i(G-2tL_{u_0})=A_t-2itT_{u_0}$$ is maximally dissipative, so that the formula of Theorem \ref{line} still holds. In particular, Theorem \ref{line} provides a formula for the extension of the flow map of the Benjamin--Ono equation to $H^s_r(\R )$ for every $s>1/2$ \cite{MP}.

\section{Final remarks}\label{final}
In the case of finite gap potentials on the torus \cite{GK}, or multisolitons on the line \cite{Sun},  formulae of Theorems \ref{torus} and \ref{line} take place in finite dimensional vector spaces, and they reduce to calculations on finite dimensional matrices, as already observed in these  references.

We expect  Theorems \ref{torus} and \ref{line} to be useful for the study of long time behaviour of solutions to the Benjamin--Ono equation. This is particular important on the line, where soliton resolution for generic data is still an open problem (see however \cite{IT} for partial results in this direction).

Let us now briefly discuss applications of a similar approach to other integrable equations. First of all, it is clear that Theorems \ref{torus} and \ref{line} easily extend to the spin Benjamin--Ono system \cite{BLL}, \cite{G}. Furthermore, these formulae could probably be very useful in the study of the small dispersion limits of these equations, in particular the half--wave maps equation \cite{GL0}, \cite{BLL}. We also expect similar formulae to hold for the recently introduced Calogero--Moser equation \cite{GL}, since the  Lax pair of operators  for this equation enjoys similar commutation properties with the shift structure of the Hardy space. Finally, as we already observed, a similar formula is known to hold for the cubic Szeg\H{o} equation on the torus \cite{GG}, and it is possible to adapt the approach with the operator $G$ developed in this paper  in order to get an explicit formula for the cubic Szeg\H{o} equation on the line \cite{GP}.
On the other hand, we have no clue whether such explicit formulae could be extended to  KdV, cubic NLS or DNLS equations.

\end{document}